\documentclass[11pt]{article}
\usepackage{amssymb,latexsym}

\newtheorem{theorem}{Theorem}[section]

\newtheorem{defi}{Definition}[section]
\newtheorem{lemma}{Lemma}[section]
\newtheorem{prop}{Proposition}[section]

\def\binom#1#2{{#1}\choose{#2}}

\def\slfrac#1#2{\hbox{\kern.1em %
 \raise.5ex\hbox{\the\scriptfont0 #1}\kern-.11em %
 /\kern-.15em\lower.25ex\hbox{\the\scriptfont0 #2}}}

\newcommand{\eqn}[1]{(\ref{#1})}

\newcommand{\eeq}{\end{equation}}
\newcommand{\beql}[1]{\begin{equation}\label{#1}}
\newcommand{\bsq}{{\vrule height .9ex width .8ex depth -.1ex }}

\newcommand{\NN}{{\mathbb N}}

\newcommand{\RR}{{\mathbb R}}
\newcommand{\TT}{{\mathbb T}}
\newcommand{\ZZ}{{\mathbb Z}}

\newcommand{\bs}{{\bf s}}

\newcommand{\sE}{{\cal E}}

\newcommand{\sP}{{\cal P}}

\newcommand{\sX}{{\cal X}}
\newcommand{\sY}{{\cal Y}}

\makeatletter
\def\@sect#1#2#3#4#5#6[#7]#8{\ifnum #2>\c@secnumdepth
     \def\@svsec{}\else
     \refstepcounter{#1}\edef\@svsec{\csname the#1\endcsname.\hskip .75em }\fi
     \@tempskipa #5\relax
      \ifdim \@tempskipa>\z@
        \begingroup #6\relax
          \@hangfrom{\hskip #3\relax\@svsec}{\interlinepenalty \@M #8\par}%
        \endgroup
       \csname #1mark\endcsname{#7}\addcontentsline
         {toc}{#1}{\ifnum #2>\c@secnumdepth \else
                      \protect\numberline{\csname the#1\endcsname}\fi
                    #7}\else
        \def\@svsechd{#6\hskip #3\@svsec #8\csname #1mark\endcsname
                      {#7}\addcontentsline
                           {toc}{#1}{\ifnum #2>\c@secnumdepth \else
                             \protect\numberline{\csname the#1\endcsname}\fi
                       #7}}\fi
     \@xsect{#5}}
\def\@begintheorem#1#2{\it \trivlist \item[\hskip \labelsep{\bf #1\ #2.}]}

\def\plain{plain}\ifx\fmtname\plain\csname fi\endcsname
     
     \input docstrip
     \preamble

     Do not distribute the stripped version of this file.
     The checksum in the header refers to the documented version.

     \endpreamble
     \generateFile{here.sty}{t}{\from{here.doc}{}}
     \endinput
\fi
\ifcat a\noexpand @\let\next\relax\else\def\next{%
    \documentstyle[here,doc]{article}\MakePercentIgnore}\fi\next
\ifx\@Hxfloat\@Hundef\else\expandafter\endinput\fi
\let\@Hxfloat\@xfloat
\def\@xfloat#1[{\@ifnextchar{H}{\@HHfloat{#1}[}{\@Hxfloat{#1}[}}
\def\@HHfloat#1[H]{%
\expandafter\let\csname end#1\endcsname\end@Hfloat
\vskip\intextsep\vbox\bgroup\def\@captype{#1}\parindent\z@
\ignorespaces}
\def\end@Hfloat{\egroup\vskip \intextsep}

\makeatother

\catcode`\@=11
\renewcommand{\section}{
        \setcounter{equation}{0}
        \@startsection {section}{1}{\z@}{-3.5ex plus -1ex minus
        -.2ex}{2.3ex plus .2ex}{\large\bf}
        }
\catcode`@=12

\begin{document}

\begin{center}
{\Large 
{\bf Benford's Law for the $3x+1$ Function}
}

\vspace{1.5\baselineskip}
{\em Jeffrey C. Lagarias} \\ 

\vspace*{1.0\baselineskip}
{\em K. Soundararajan } \\

\vspace*{.5\baselineskip}
Dept. of Mathematics \\
University of Michigan \\
Ann Arbor, MI 48109-1109\\
\vspace*{1.5\baselineskip}

(August 9, 2006) \\
\vspace{3\baselineskip}
{\bf ABSTRACT}
\end{center}
Benford's law (to base B) for an infinite sequence $\{x_k: k \ge 1\}$
of positive quantities $x_k$ 
is the assertion that $\{ \log_B x_k : k \ge 1\}$ is
uniformly distributed $(\bmod~ 1)$.
The $3x+1$ function $T(n)$ is given by 
$T(n)=\frac{3n+1}{2}$ if $n$ is odd, and $T(n)= \frac{n}{2}$
if $n$ is even.
This paper studies the initial iterates 
$x_k= T^{(k)}(x_0)$ for $1 \le k \le N$ of the $3x+1$
function, where $N$ is fixed. It shows that for most 
initial values $x_0$, such sequences 
approximately satisfy Benford's law, in the sense that the discrepancy of
the finite sequence $\{ \log_B x_k: 1 \le k \le N \}$ is small. \\


%***************************************************
%
%
% Section 1 Introduction
%
%
%***************************************************
\setlength{\baselineskip}{1.0\baselineskip}

\section{Introduction}

The $3x+1$ problem concerns the behavior under iteration of the
map $T: \ZZ \to \ZZ$ given by $T(n) = \frac{n}{2}$
or $T(n) = \frac{3n+1}{2}$ according as $n$ is even or odd. 
That is, $T(2m)=m$ and $T(2m+1)=3m+2$. The notorious $3x+1$ Conjecture
asserts that when started from any positive integer $n$,
some iterate $T^{(k)}(n)=1$; it remains unsolved. 
Surveys of work on this problem appear in Lagarias \cite{La85}
and Wirsching \cite{Wi98}.

It is well known that the initial iterates of this map exhibit
a ``random'' character. This holds in the sense that the
initial  iterates of
a randomly selected integer appear to be even or odd
with equal probability. Such a result can be rigorously
justified if one takes the interval $1 \le n \le X= 2^k$
and considers only the first $k = \log_2 X$ iterations
(see \cite[Theorem A]{La85}). This leads to the
rapid decay of most trajectories of the iteration under
$T$, at an exponential rate, with an expected decrease
by a multiplicative factor 
$\sqrt{\frac{3}{4}} \approx 0.86602$      
at each step. 
These facts support the conjecture that all orbits of the $3x+1$
iteration enter a bounded set, and hence fall into a finite
number of periodic orbits. Heuristic stochastic models
(Lagarias and Weiss \cite{LW92}, Borovkov and Pfeifer \cite{BP00})
predict that for an integer of size about $X$
the ``random'' character above persists for about the first 
$\alpha \log X$ iterates, with 
$\alpha=  \left(\frac{1}{2} \log \frac{3}{4}\right)^{-1}\approx 6.95212$;
the model predicts  
most integers of size near $X$ will
arrive at the periodic orbit $\{ 1, 2\}$ near this number of iterations.
The stochasic model in \cite{LW92} also predicts that
for large $n$ the number of steps to enter a periodic orbit should never
exceed $42 \log n$.  Experimentally, Roosendaal \cite{Ro04} 
has found a number $n$ of size $7.2 \cdot 10^{21}$ which 
requires about $36.7 \log n$ iterations before entering the 
periodic orbit $\{1,2\}$.

The present paper concerns the base $B$ expansion of
the initial sequence of
the first $N$ iterates of the
$3x+1$ map on a random starting value $n$, drawn
from $1 \le n \le X$ where $X \ge 2^N$.
This  is in the region of the dynamics
 where most trajectories
are decreasing at an exponential rate, before they
enter a periodic orbit.  It
shows that, in a certain sense,
 the leading digits of the base $B$ expansion of
most such sequences approximately satisfy 
a strong form of Benford's law.
Here
Benford's law concerns the distribution of the initial digits in the
base $B$ expansion of an infinite sequence 
$\sX=\{x_{1}, x_{2}, x_{3},...\}$ of positive real numbers.
The original version of Benford \cite{Be38} in 1938 concerned the
first few leading digits in the decimal expansion of real numbers in tables;
the distribution had already
been formulated by Newcomb \cite{New81} in 1881.
 An infinite sequence $\sX$ is said to satisfy the
{\it strong  Benford's law} (to base $B$) if for each fixed $k \ge 1$,
the first $k$ digits in the $B$-ary
expansion of $\{x_{1}, x_{2}, x_{3},...\}$ approach limiting
probabilities given by the ``$B$-ary Benford distribution'',
which we specify below. This is known to be  equivalent
to the condition
that the  associated infinite sequence
$ y_i := \log_B x_i$ is uniformly distributed modulo one
(Diaconis \cite[Theorem 1]{Di77}). In what follows
we adopt this criterion as our
definition of Benford's law.

This paper is motivated by
work of Kontorovich and Miller \cite{KM04}, who  showed 
that certain statistics drawn from $3x+1$ iterates 
approximately obey Benford's law. They treated a version
of the $3x+1$ iteration in which the initial starting
point  $w_0$ is an odd integer, 
and they studied the subset of  the successive odd integers 
$\{w_1, w_2, ...\}$ appearing in the $3x+1$ iteration of $w_0$.
Here $w_i = T^{(k_i)}(w_0)$ where $k=k_i$
is the $i$-th  value  where
$T^{(k)}(w_0)$ is odd. They showed that for a suitable
natural initial distribution on the  odd integers 
drawn from $1 \le w_0 \le X$,
and for a suitable  number $k$ of iterates (growing slowly with $X$),
as $X \to \infty$ the distribution of the $B$-ary
digits of the ratios $w_k/w_0$ 
approached the $B$-ary Benford distribution,
 provided that $B$ was not a power of $2$.
More precisely, they obtained the Benford  distribution in a double limit,
in which $X \to \infty$ with $k$ held fixed, and after this
taking $k \to \infty$. They also gave results of numerical simulations
indicating that 
the distribution of the odd $3x+1$ iterates $\{w_1, w_2, ..., w_k\}$
starting from an odd $w_0$
themselves should approximately
satisfy Benford's law,
for all integer bases $B$ not a power of $2$.
In the case where $B$ is a power of $2$, they showed
that a  double limiting
distribution  exists,  
but is  not the $B$-ary Benford distribution.

The main result of this paper, Theorem~\ref{th21}  in \S2,  
establishes in
a quantitative form the assertion
that most initial sequences of  the first $N$ iterates of the $3x+1$ 
function  approximately satisfy the strong Benford
law. It applies to a finite sequences of initial
$3x+1$ iterates $\{x_1, x_2, ..., x_N\}$,
and obtains an upper bound on 
the discrepancy  $D( \{y_1, y_2, ..., y_N\})$  of the sequence of numbers
 $y_j= \log_{B} x_j$ for most such sequences.
The discrepancy  is a well-known statistic which is a 
measure of distance to the uniform distribution.
It is defined in \S2, and relevant  properties of 
discrepancy are treated in \S3.  We obtain an explicit 
upper bound on the number
of ``exceptional'' sequences for 
which the discrepancy is large.
We treat  $3x+1$ iterates including both even
and odd iterates, and our main result implies  
convergence to a generalized
Benford's law for all bases $B \ge 2$, including $B$ being a
power of $2$.
The anomalous behavior of powers of 
$2$ in the results of Kontorovich and Miller \cite{KM04}
is associated to their restriction to the subset of iterates that are
odd integers. 

The basic approach is as follows. We use the fact that the
initial iterates of a large randomly chosen integer $n$ are
well approximated by a stochastic process that takes
 $T(n) = \frac{n}{2}$ or $T(n) = \frac{3n}{2}$ with
equal probability. Taking logarithms to the base $B$, we are
reduced to studying the stochastic process
which sets either  
$$
y_{n+1} = y_{n} + \theta_1,
$$
or
$$
y_{n+1} = y_{n} + \theta_2
$$
with equal probability, where 
$$
\theta_1= \log_B \frac{3}{2} ~~\mbox{and}~~ \theta_2= \log_B \frac{1}{2}.
$$
In \S4 we consider this process in its own right, for arbitrary 
$(\theta_1, \theta_2)$.
We first show  that the  realizations
$$\omega= \{ y_{n}: n= 1,2,3,...\}$$
 of such a stochastic process
for general $(\theta_1, \theta_2)$ 
are uniformly distributed modulo one with probability one,
if and only if at least one of 
$\theta_1$ or $\theta_2$
is irrational. 
The main result of \S4 shows that
if  the numbers $\theta_1$
and $\theta_2$ are not 
simultaneously  well approximable by rational 
numbers, as specified  by a two-dimensional ``Diophantine property'',
then for any fixed $N$ most initial segments of
length $N$ are close to the uniform distribution, quantitatively
given by an upper bound on their discrepancy.

In \S5 we apply the results of \S4 to the $3x+1$ iteration.
We show using a result of Rhin \cite{Rh87} 
that $\theta_1= \log_B \frac{3}{2}$
and $\theta_2= \log_B \frac{1}{2}$ have suitable 
two-dimensional Diophantine
properties for the results in \S4 to apply. Then we
establish  that the $3x+1$ iterates are sufficiently close
to realizations of the stochastic process to obtain
upper bounds on the discrepancy of sequences for
most initial inputs, provided we average over $1\le n \le X$,
and for $N$ iterates we require $X \ge 2^{N}$.
Putting all these results
 together yields  the main result, Theorem~\ref{th21}. 
 
 The main result is established here  for the
$3X+1$ function, but the methods used apply equally well
to number-theoretic
maps of a similar nature, such as the $Qx+1$ function, for
odd $Q$, with
$T_{Q}(n)= \frac{n}{2}$ or $T_{Q}(n)= \frac{Qn+1}{2}$
according as $n$ is even or odd.   Results analogous to Theorem~\ref{th21}
should hold for the distribution of the first $N$ iterates of such functions.
 For $Q \ge 5$ it is expected that most initial values 
of the  $Qx+1$ iteration never enter a periodic orbit, but diverge 
to $+\infty$.
It seems possible that the infinite sequence $\{ x_{n}: n \ge 0\}$
of a divergent orbit might actually satisfy a strong Benford's law. However
at present there seems no approach to  address this question; 
even the existence of
a divergent orbit for the $Qx+1$ function, for any $Q \ge 5$, remains an
open problem.

There has been other work showing that 
the iterates of certain dynamical systems
satisfy Benford's law, see Berger, Bunimovich
and Hill \cite{BBH05} and Berger \cite{Be05}.  For various properties 
of Benford's law, see Hill ~\cite{Hi95}, \cite{Hi96}.
Finally we observe that the 
approach of Kontorovich and Miller \cite{KM04}
to Benford's law for $3x+1$ iterates introduced several
ideas to this problem,
including approximation to a stochastic process
(not the one studied here), as well as a relation
to Diophantine properties of certain constants.
Their approach starts
 from a structure formula for odd iterates of the  $3x+1$
function given by Sinai \cite{Si03}
and extended in Kontorovich and Sinai \cite{KS03}
to a wider class of maps. 
Their main result  (\cite[Theorem 5.3]{KM04}) for the
$3x+1$ function
establishes the uniform
distribution in a double limit 
of $y_i := \log_{B} (\frac{w_i}{w_0})$ for any
real base $B$ such that $\log_B 2$ satisfies
a one-dimensional
Diophantine property, as defined
in \S4 below.

\paragraph{Notation.} We let $\lfloor x \rfloor$ denote
the largest integer that does not exceed $x$, and we
let $\{\{x \}\} := x - \lfloor x \rfloor$ denote the
fractional part of $x$, with $0 \le \{\{x \}\} < 1$. 
Finally $\Vert x\Vert = \min_{n\in {\Bbb Z}} |n-x|$ 
denotes the distance of $x$ from its nearest integer. 

%***************************************************
%
%
% Section 2 Main Result
%
%
%***************************************************

\section{Main Result}

Benford's law concerns the distribution of the initial digits in the
base $B$ expansion of an infinite sequence 
$\sX=\{x_{1}, x_{2}, x_{3},...\}$ of positive real numbers.
An infinite sequence is said to satisfy the
{\em strong Benford's law (to base B)} 
if the associated infinite sequence $\sY= \{ y_1, y_2, y_3,...\}$ 
given by the base $B$ logarithms 
$ y_i := \log_B x_i$ is uniformly distributed modulo one.
Suppose that the numbers $x_n$ have $B$-ary expansion
$$
x_n= B^{M_n} \left(\sum_{k=0}^{\infty} d_k^{(n)} B^{-k}\right)
$$
with $1 \le d_0^{(n)} \le B-1$ and $0 \le d_k^{(n)} \le B-1$ for $k \ge 1$.
Benford's law is the statement that
$$
{\rm{Prob}}\left[ d_0^{(n)} = d \right] 
= \log_B (d+1) - \log_{B} d
$$
for $1 \le d \le B-1$, in which the ``probability'' is interpreted
as a limiting frequency in the first $N$ values of $x_n$ as $N \to \infty$. 
More generally the  {\em strong Benford
probability} of observing a given block of $K$
digits $[d_0 d_1 ...d_{K-1}]$, with $d_0 \ne 0$,
is given by 
$$
{\rm{Prob}}\left[ d_0^{(n)}d_1^{(n)}\cdots d_{K-1}^{(n)} :=
d_0 d_1\cdots d_{K-1} \right]
= \log_B (r+B^{-K+1}) - \log_{B} r,
$$
where
\beql{228a}
r = \sum_{j=0}^{K-1} d_j B^{-j}.
\eeq

The departure from uniform distribution modulo one of
a finite set $\sY$ can be measured using the discrepancy.

\begin{defi}~\label{de21}
The discrepancy $D(\sY)$ of  
a finite set $\sY =\{ y_1, y_2, ..., y_N\}$ of real numbers is
defined as follows. For $0 \le \alpha \le \beta \le 1$ set
\beql{240} 
Z(\sY; \alpha, \beta) := \frac{1}{N}\# \{i : \alpha  \le \{\{y_i\}\} 
\le \beta\}.
\eeq
in which $\{\{ y \}\}= y - \lfloor y \rfloor$ is the fractional part of $y$,
and then let
\beql{241}
D(\sY; \alpha, \beta) := Z(\sY; \alpha, \beta) - 
\left( \beta - \alpha \right) .
\eeq
The  {\it (normalized) discrepancy} $D(\sY)$ is  then
\beql{242} 
D(\sY) := \sup_{0 \le \alpha \le \beta \le 1} |D(\sY; \alpha, \beta)|.
\eeq
It is also given by
\beql{242b}
D(\sY) = 
\sup_{0 \le \alpha \le 1} D(\sY; 0, \alpha)
- \inf_{0 \le \alpha \le 1} D(\sY; 0, \alpha).
\eeq
\end{defi}

One has $0 \le D(\sY) \le 1$;
smaller values of  $D(\sY)$ correspond to more  uniformly
spaced sets  $\sY$  modulo one.
No finite distribution can be perfectly uniform, so there
is a nonzero lower bound on the discrepancy of all sequences of 
length $N$. This minimal value of the discrepancy
is attained by equally spaced elements $y_i = \frac{i}{N}$
for $0 \le i \le N-1$, with $D(\sY) = \frac{1}{N}.$
This notion of discrepancy is {\it  translation-invariant}; that is,
for any real $y_0$, one has  
\beql{242bb}
D(\sY + y_{0}) = D(\sY).
\eeq
Some authors treat instead a
(normalized)  non-translation invariant discrepancy
$$
D^{\ast}(\sY) :=
 \sup_{0 \le \alpha \le 1} \Big| Z(\sY; 0, \alpha) - \alpha\Big|.
$$
This is related to $D(\sY)$ by the inequlities
$D^{\ast}(\sY) \le D(\sY) \le 2 D^{\ast}(\sY)$.

Our  definition of discrepancy follows 
Kuipers and Niederreiter \cite{KN74} and 
Drmota and Tichy \cite{DT97}).
A few authors (Montgomery \cite{Mo94})
study an unnormalized discrepancy that does not divide
by $N$; this version of the discrepancy takes values
between $0$ and $N$.

The main result of this paper is an upper bound on discrepancy
of the base $B$ logarithms of most initial $3x+1$ sequences.

%%********************************************
%
% theorem 2.1 statement
%
%********************************************
\begin{theorem}~\label{th21}
Let $B \ge 2$ be a fixed integer base. 
 For each $N \ge 1$ 
and each $X \ge 2^N$, most initial seeds $x_0$ in 
$1 \le x_0 \le X$ have first $N$ initial 
 $3x+1$ iterates $\{ x_k: 1 \le k \le N\}$
that satisfy the discrepancy bound 
\begin{equation}~\label{eq101}
D( \{ \log_B x_k: 1 \le k \le N \} ) \le 2 N^{ - \frac{1}{36}}.
\end{equation}
The set $\sE(X, B)$
of exceptional initial seeds $x_0$ in $1 \le x_0 \le X$ 
that do not satisfy this bound has cardinality 
\beql{101b}
|\sE(X, B)|  \le c(B) {N^{-\frac{1}{36}}}X,
\eeq
where $c(B)$ is a positive constant depending only on $B$.
\end{theorem}

This result implies approximation to base $B$ Benford's 
law, as follows. Let $\sX = \{ x_{1}, ..., x_{N}\}$ be a
set of positive real numbers, and set $y_{i}= \log_{B} x_{i}$
and $\sY= \{ y_{1}, ..., y_{N}\}$. Let $1\le r< B$ be a $B$-ary rational 
as in \eqn{228a} with $1 \le r < B$.  Requiring that the 
 first $K$ digits of $x_n$ match the digits of $r$ is clearly equivalent to 
having $\{\{ y_n\}\}$ lie in the interval
 $[\log_B r, \log_B (r+B^{-K+1}))$.  
Fron the definition of discrepancy, we have that 
$$
\Big|~\frac{1}{N}\# \left\{ 1 \le i \le N:~ \log_{B} r \le 
\{\{ \log_B x_i\}\} < \log_B(r+B^{-K+1})  \right\} - 
\log_{B} \Big(\frac{r+B^{-K+1}}{r}\Big) \Big|
$$
is bounded above  by $D(\{ y_1, y_2, ..., y_N \})$,
independent of $K$. Theorem~\ref{th21}
upper bounds this discrepancy for the intial iterates of
 most $3x+1$ sequences.

%***************************************************
%
%
% Section 3 preliminaries on discrepancy/exponential sums
%
%
%***************************************************

\section{Discrepancy and Exponential Sums}

We will use standard criteria for uniform
distribution of an infinite sequence $\sY =\{y_1, y_2, ...\}$ 
in terms of exponential sums
and of the discrepancy of its initial segments (\cite[Chap. 1]{Mo94}).

For an infinite sequence $\sY =\{y_1, y_2, ...\}$ 
we let $\sY_N$ denote the first $N$
elements of $\sY$.  For integers $k$, we associate to $\sY_N$ the
`Fourier coefficients'
\beql{300}
\hat{U}_N (k, \sY) = \hat{U} (k, \sY_N) := \sum_{j=1}^N e^{2 \pi i k y_j}.
\eeq

\begin{prop}~\label{pr31}
For an infinite
sequence $ \sY = \{ y_1, y_2, ... \}$ of real numbers,
the following conditions on $\sY$ are equivalent.

(1) The sequence $\sY$ is uniformly distributed modulo one.

(2) (Weyl's criterion) For each nonzero integer $k$ we have 
\beql{301}
\lim_{N \to \infty} \frac{1}{N} |\hat{U}_N(k, \sY)| = 0.
\eeq

(3) For any  properly Riemann integrable function $F$ on $[0,1]$,
\beql{302} 
\lim_{N \to \infty} \frac{1}{N} \sum_{j=1}^N F(y_j) = \int_{0}^1 F(t)dt.
\eeq

(4) The discrepancy $D(\sY_N)$ satisfies  
\beql{303}
\lim_{N \to \infty}  D \left( \{y_1, y_2 ,..., y_N \}\right) = 0.
\eeq
\end{prop}

\paragraph{Proof.} Here (1)-(3) are Weyl's criterion in
\cite[page 1]{Mo94}, and the equivalence of (1) and (4) appears 
in \cite[page 2]{Mo94}.
$~~~\bsq$ \\

We will need a quantitative relation between exponential sums
$\hat{U}_N(k, \sY)$ and discrepancy, given by the Erd\H{o}s-Turan
inequality.

\begin{prop}~\label{pr32}
(Erd\H{o}s-Turan Inequality)
For any positive integer $K \ge 1$,
\beql{310}
D \left( \{y_1, y_2 ,..., y_N\}\right) \le \frac{1}{K+1}+ 3 \sum_{k=1}^K \frac{1}{k} \Big|\frac{1}{N}\sum_{n=1}^N e^{2 \pi i k y_n}\Big|
\eeq
\end{prop}

\paragraph{Proof.} This is a weak form
of the Erd\H{o}s-Turan inequality.
A short proof of it is  given  in Montgomery \cite[page 8]{Mo94}
(after normalizing the discrepancy). For a stronger form see
 Kuipers and Neiderreiter \cite[Ch. 2, Theorem 2.5]{KN74}.
$~~~\bsq$ \\

We will also need 
the following simple bound on the change in discrepancy under perturbation.

\begin{prop}~\label{pr33}
If $|y_{i}- \tilde{y}_{i}| \le \epsilon$ for $1 \le i \le N$ then
\beql{331}
| D\left( \{y_1, y_2 ,..., y_N\}\right) -
D\left( \{\tilde{y}_1, \tilde{y}_2 ,..., \tilde{y}_N\}\right)| 
\le 2 \epsilon.
\eeq
\end{prop}
\paragraph{Proof.} 
Let $\sY$ and $\sY'$ denote the sets in the Proposition. Suppose first the discrepancy 
$D(\sY)$ is attained on an interval $J=[\alpha, \beta]$ with $Z(\sY;J) -|J|>0$.  
If $\alpha >\epsilon$ and $\beta <1-\epsilon$, then with $J^{\prime}=[\alpha-\epsilon,
\beta+\epsilon]$ we see that $Z(\sY';J') \ge Z(\sY;J)$ and it 
follows that 
\[ 
D(\sY')\ge Z(\sY'; J') - |J'| \ge Z(\sY; J) - |J| - 2\epsilon = D(\sY)-2\epsilon. 
\]
 If $\alpha <\epsilon$ or $\beta>1-\epsilon$ 
we would still like to consider $J^{\prime}\subset [0,1]$ which is 
the image $\pmod 1$ of the interval $[\alpha-\epsilon,\beta+\epsilon]$. 
The only issue is that $J'$ now consists of two intervals, one near $0$ and the other near $1$.  
However, the complement $J'^c$ is a genuine interval and we have 
$|J'^c| - Z(\sY';J'^c)= Z(\sY';J')-|J'|\ge D(\sY)-2\epsilon$.  Thus 
we have again that $D(\sY') \ge D(\sY)-2\epsilon$.  

In the remaining case that the discrepancy $D(\sY)$ is attained on an interval $J=[\alpha,\beta]$ 
with $|J|-Z(\sY;J)>0$, we consider $J'=[\alpha+\epsilon,\beta-\epsilon]$ if $\beta-\alpha>2\epsilon$, and 
$J'$ to be the empty interval otherwise.  We deduce in this case also that $D(\sY')\ge D(\sY)-2\epsilon$.  

Since $\sY$ and $\sY'$ are interchangeable
in the argument, 
we obtain $D(\sY) \ge D(\sY')- 2 \epsilon$, completing the proof.  $~~~\bsq$ \\

In the sequel we will obtain bounds on exponential sums and
from this derive bounds on the discrepancy using the Erd\H{o}s-Turan
inequality. We will approximate the values $y_i =\log_B x_i$
of the $3x+1$ iterates of a randomly drawn initial value $x_0$
by the values of a stochastic process, of a type which
we analyze in the next section.

%***************************************************
%
%
% Section 4 Stochastic Process
%
%
%***************************************************

\section{Stochastic Process}

We study the following family of stochastic processes.
We suppose that we are given two real numbers $(\theta_{1}, \theta_{2})$,
and an initial value $y_0$. 
The discrete stochastic process 
$\sP(\theta_1,\theta_2, y_{0})$ has
 realizations of the form
\beql{400c}
\omega=( y_{1}, y_{2}, y_{3},...)
\eeq
in which the $y_{i}$ are generated
from the initial value $y_{0}$ by  choosing  
\beql{401}
y_{n+1} = 
y_{n}+ \theta_{1} {\rm {\ with \ probability \ }} \frac{1}{2}, \ \ 
{\rm and \ \ } y_{n+1}= 
y_{n}+ \theta_{2}  {\rm {\ with \ probability \ }} \frac{1}{2},
\eeq
where each step is an independent Bernoulli trial. 
We think of the $y_{i}$ as given modulo one, in which case 
this process is a Bernoulli mixture of two rotations of the circle.

%*******************************************
% Theorem 4.1
%*******************************************
%
\begin{theorem}~\label{th41}
If at least one of  $\theta_1$ or $\theta_2$ is irrational, 
then for any fixed initial value $y_{0}$
the process $\sP(\theta_1, \theta_2, y_{0})$ has a probability one subset of 
realizations $\omega=(y_{1}, y_{2}, ...)$  that are
uniformly distributed modulo one. Equivalently, with probability one,
\beql{403}
\lim_{N \to \infty}  D \left( \{ y_{1},...,y_{N} \} \right)=0.
\eeq
\end{theorem}

Note that if $\theta_1$ and $\theta_2$ are both rational numbers,
then the values $y_i$
can only take  a finite number of distinct values 
modulo one and no realization $\omega$ is uniformly distributed modulo one. 
We also remark that Theorem 4.1 may be easily generalized to cover 
Bernoulli mixtures of $K$ rotations of the circle.  

Theorem~\ref{th41} will be derived using exponential sums. 
We first study finite initial segments of length $N$ of such a stochastic
process $\sP(\theta_1, \theta_2,y_{0})$.
We let
$$
\omega_{N}:=(y_{1},y_{2},...,y_{N}).
$$
denote such an initial segment, and write 
${\Bbb E}_{{\omega_{N}}}[f(\omega_{N})]$
for the expected value of a random variable over 
the process restricted to  these initial segments.
We begin by calculating the second moment of 
the individual Fourier coefficients $\hat{U}_N(k, \omega)$ of 
$\omega_N$.

%*******************************************
% Lemma 4.1
%*******************************************
%

\begin{lemma}~\label{le41}
For each $N \ge 1$ and each $k \in {\Bbb Z}$ 
\beql{420}
{\Bbb E}_{\omega_N} \Big[ |\hat{U}_N(k, \omega)|^2\Big]
= N +  2 {\rm Re \ } \Big( \sum_{r=1}^{N} (N-r) 
\Big(\frac{e^{2\pi ik \theta_1}+e^{2\pi i k\theta_2}}{2}\Big)^r\Big).
\eeq
If at least one of $\theta_1$ or $\theta_2$ is irrational, 
then for each non-zero integer $k$ and each $N\ge 1$ 
\beql{431}
{\Bbb E}_{{\omega_{N}}} \Big[ |\hat{U}_{N} (k, \omega) |^{2} \Big] \le
\Big(1+ \frac{8}{|2-e^{2\pi ik\theta_1}-e^{2\pi ik\theta_2}|}\Big) N 
\le  \Big( 1 + \frac{1}{\Vert k\theta_1\Vert^2 +
\Vert k\theta_2\Vert^2}\Big) N,
\eeq
where $\Vert \xi \Vert = \min_{n\in {\Bbb Z}} |\xi-n|$ denotes the 
distance between $\xi$ and its nearest integer.  
\end{lemma}

\paragraph{Proof.} Observe that 
$$
|{\hat U}_N(k,\omega)|^2 = \Big| \sum_{j=1}^{N} e^{2 \pi ik y_j} \Big|^2 
= N + 2 \ {\rm Re \ } \Big(\sum_{1\le j < \ell \le N} e^{2\pi i k 
(y_\ell-y_{j})}\Big).
$$
If we write $r=\ell-j$ then $y_\ell -y_j$ is a 
sum of $r$ random variables each taking the values $\theta_1$ or 
$\theta_2$ with equal probability.  Thus 
$$
{\Bbb E}_{\omega_N}\Big[ e^{2\pi i k(y_\ell - y_j)} \Big] 
= \Big( \frac{e^{2\pi ik\theta_1}+e^{2\pi i k\theta_2}}{2} 
\Big)^{\ell -j}.
$$
Since for $1\le r\le N$ there are $N-r$ pairs $1\le j <\ell \le N$ 
with $\ell-j=r$, we conclude that 
$$
{\Bbb E}_{\omega_N}\Big[|{\hat U}_N(k,\omega)|^2 \Big] 
=  N + 2 \ {\rm Re \ }\Big( \sum_{r=1}^{N} (N-r) \Big(
\frac{e^{2\pi ik\theta_1}+e^{2\pi i k\theta_2}}{2} \Big)^r\Big).
$$
This proves \eqn{420}.  

For any $z\neq 1$ we note that 
$$
\sum_{r=1}^{N} (N-r) z^{r} = \frac{(N-1)z-Nz^2+z^{N+1}}{(1-z)^2},
$$
and so, if $|z|\le 1$ and $z\neq 1$ we get that 
\beql{4201}
\Big| \sum_{r=1}^{N} (N-r) z^{r}\Big| \le 
\frac{N|z-z^2|+|z-z^{N+1}|}{|1-z|^2}\le 
\frac{2N}{|1-z|}.
\eeq
If at least one of $\theta_1$ or $\theta_2$ is 
irrational, then for non-zero $k$ we have that $e^{2\pi ik\theta_1}+e^{2\pi 
ik\theta_2} \neq 2$, and of course $|e^{2\pi i\theta_1}+e^{2\pi i\theta_2}|
\le 2$.  Combining \eqn{420}, and \eqn{4201} with $z=(e^{2\pi ik\theta_1}
+e^{2\pi ik\theta_2})/2$, we obtain that 
$$
{\Bbb E}_{\omega_N}\Big[ |{\hat U}_N(k,\omega)|^2\Big] \le 
\Big( 1 +\frac{8}{|2-e^{2\pi ik\theta_1} -e^{2\pi i k\theta_2}|}\Big) N.
$$
For $|\xi| \le 1/2$ note that $\sin^2 (\pi \xi) \ge 4\xi^2$ and so 
\begin{eqnarray*}
|2-e^{2\pi ik\theta_1} -e^{2\pi i k\theta_2}| &\ge &
2 - \cos (2\pi  k \theta_{1}) - \cos (2\pi k \theta_{2}) \\
%& = & 2 - \cos^{2}( \pi k \theta_{1})- \cos ^{2}(\pi k \theta_{2})
%+ \sin ^{2}(\pi k \theta_{1}) + \sin ^{2}(\pi k \theta_{2}) \\
& = & 2\left(\sin^2 (\pi k \theta_1) + \sin^2 (\pi k \theta_2) \right) 
\ge 8 \left(\Vert k\theta_1 \Vert^2 + \Vert k\theta_2 \Vert^2 \right),
\end{eqnarray*}
which completes the proof of \eqn{431}.  $~~~\bsq$

%******************************************
% Proof of theorem 4.1
%*******************************************
%
\paragraph{Proof of Theorem~\ref{th41}}
We suppose that at least one of $\theta_1$ or $\theta_2$ is irrational.
We claim that for each nonzero $k$ there holds
\beql{493}
{\rm Prob}_{\omega}\Big[ \lim_{N \to \infty} 
\frac{1}{N}|\hat{U}_{N}(k, \omega)|=0\Big]
=1.
\eeq
Thus, for each fixed non-zero integer $k$,
there is a probability one set of $\omega$ such that
 
\noindent $\lim_{N \to \infty} \frac{1}{N} |{\hat U}_N(k,\omega)| =0$.  
Since the set of non-zero integers $k$ is countable, it 
follows that the set of all $\omega$ for which 
$\lim_{N \to \infty} \frac{1}{N} |{\hat U}_N(k,\omega)|=0$ 
holds simultaneously for all non-zero integers $k$ still has 
probability one. (Its complement is a countable union of
sets of measure zero.)
Now by Weyl's criterion (Proposition ~\ref{pr31}(2)) all such $\omega$
are uniformly distributed modulo one.  Proposition ~\ref{pr31}(4) then
yields  \eqn{403} with probability one.

To prove \eqn{493} it suffices to show that for each $1\ge \delta > 0$
\beql{494}
P_{\delta} := {\rm Prob}_{{\omega}}\Big[ \limsup_{N \to \infty}
\frac{1}{N}| \hat{U}_{N}(k, \omega) | \ge \delta \Big]  = 0.
\eeq
For $j\ge 1$  set 
$N_j :=\left\lfloor \frac{1}{(1- \frac{\delta}{2})^{j}} \right\rfloor$.  If 
$N_j \le N < N_{j+1}$ is such that $|{\hat U}_N(k,\omega)| 
\ge \delta N$, then we see that 
$$
|{\hat U}_{N_j}(k,\omega)| 
\ge |{\hat U}_N(k,\omega)| - \Big|
\sum_{\ell =N_{j}+1}^N e^{2\pi i k y_{\ell}}\Big| 
\ge \delta N - (N-N_j) \ge N_{j} - 
\Big(\frac{1-\delta}{1-\frac{\delta}{2}}\Big) N_{j} \ge
\frac{\delta}{2} N_j.
$$
Therefore, for any $B \ge 1$,  
\beql{495}
P_{\delta} \le {\rm Prob}_{\omega} \Big[ \limsup_{j\to 
\infty} \frac{1}{N_j} |{\hat U}_{N_j}(k,\omega)| \ge 
\frac{\delta}{2}\Big] 
\le \sum_{j=B}^{\infty} {\rm Prob}_{\omega}\Big[ 
|{\hat U}_{N_j}(k,\omega)| \ge \frac{\delta N_j}{2}\Big].
\eeq

Now 
$$ 
{\rm Prob}_{\omega}\Big[ 
|{\hat U}_{N_j}(k,\omega)| \ge \frac{\delta N_j}{2}\Big] 
\le \Big( \frac{\delta N_j}{2}\Big)^{-2} 
E_{\omega} \Big[ |{\hat U}_{N_j}(k,\omega)|^2 \Big],
$$
and by Lemma \ref{le41} this is 
$$
\le \frac{4}{\delta^2} \Big(1 +\frac{1}{\Vert k\theta_1\Vert^2 +
 \Vert k\theta_2 \Vert^2}\Big) \frac{1}{N_j}.
$$
We use this in \eqn{495}, and obtain that 
for any $B\ge 1$, 
$$
P_{\delta} \le \frac{4}{\delta^2}
 \Big(1 +\frac{1}{\Vert k\theta_1\Vert^2 + \Vert k\theta_2 \Vert^2 }\Big) 
  \sum_{j=B}^{\infty} \frac{1}{N_j}.
$$
Since the $N_j$ grow exponentially, letting $B\to \infty$ 
we may conclude that $P_\delta=0$.  This establishes \eqn{494}, 
and \eqn{493} and the Theorem follows.
$~~~\bsq$ \\

For general non-rational pairs $(\theta_1,\theta_2)$ 
the convergence rate to zero in
\eqn{403}, or equivalently \eqn{493}, can be arbitrarily slow. 
To obtain explicit
bounds on the convergence rate in \eqn{403} one must impose restrictions
 on the Diophantine approximation properties 
of the numbers $\theta_1$ and $ \theta_2$.  
The following definition has been much used in connection with
``small divisors''
problems in dynamical systems, 
cf. Herman \cite{He79}, Yoccoz \cite{Yo84}, \cite{Yo02},
and  in number theoretical dynamics cf. Marklof \cite{Ma03}.

%*******************************************
% Defn 4.1
%*******************************************
%

\begin{defi}~\label{de41}  A real number $\theta$ is said to be 
{\rm Diophantine with exponent $\alpha$}
if there is a positive constant $C(\theta)$ such that
for all integers $k\ge 1$ 
\beql{404b}
 \Vert k\theta \Vert  \ge C(\theta) |k|^{- \alpha}.
\eeq
\end{defi}

Any real number that is Diophantine with some positive  exponent $\alpha$
is  irrational; necessarily $\alpha \ge 1$.
For any $\alpha > 1$, the set of real numbers that are Diophantine with
exponent $\alpha$ has full Lebesgue measure. In fact the exceptional set of
real numbers that are not Diophantine with a given exponent $\alpha >1$ has
Hausdorff dimension $f(\alpha)$ with $f(\alpha) < 1$. {\it Liouville numbers}
are those real numbers that are not Diophantine for any finite exponent, and
they form an uncountable set of Hausdorff dimension zero.  The
set of real numbers that are Diophantine with exponent $\alpha=1$
comprise the {\it badly approximable numbers}, and these form a set
of Hausdorff dimension one but Lebesgue measure zero. 

 In this paper we use the following generalization of 
this notion to simultaneous
 approximation, which is the complement of the notion of
$d$-dimensional very well approximable vectors  appearing
in Schmidt \cite{Sch71}.
  
%*******************************************
% Defn 4.2
%*******************************************
%

\begin{defi}~\label{de42}  The vector $(\theta_1,\theta_2,.., \theta_{d})$
of real numbers is said to be 
{\rm $d$-dimensional Diophantine with exponent $\alpha$}
if there is a positive constant $C(\theta_1,\theta_2,...,\theta_{d})$ 
such that
for all integers $k\ge 1$ 
\beql{404c}
\max ( \Vert  k\theta_1 \Vert, \Vert k \theta_2 \Vert,
..., \Vert k\theta_{d}\Vert)
  \ge C(\theta_1,\theta_2,...,\theta_{d}) k^{-\alpha}.
\eeq
\end{defi}

This notion has been used  in
the dynamical system context by Marklof \cite{Ma02}, \cite{Ma05}.
Here we use the case $d=2$.
The multidimensional notion is less restrictive
than the case $d=1$ in the sense that if any $\theta_{i}$
is one-dimensional Diophantine with exponent $\alpha$,
then the vector $(\theta_{1}, ..., \theta_{d})$ will
be $d$-dimensional Diophantine with the same or smaller exponent.

%*******************************************
% Theorem 4.4
%*******************************************
%
 
The next result gives bounds on the expected size of
the discrepancy of a finite initial segment of this stochastic process,
under suitable Diophantine conditions on $(\theta_1,\theta_2)$.

%**************************************************************
%
%Theorem 4.2
%
%**************************************************************
\begin{theorem}~\label{th42}  Suppose that the 
pair  $(\theta_1,\theta_2)$ is two-dimensional
Diophantine with exponent $\alpha$.  
Then there is  a constant $C_{2}(\theta_1,\theta_2)$ such that
for all $N \ge 1$,
\beql{451}
{\Bbb E}_{\omega_{N}}[  D(\{y_{1}, y_{2}, ..., y_{N}\})  ]\le
C_{2}(\theta_1,\theta_2) N^{- \frac{1}{2(1+\alpha)}}.
\eeq
\end{theorem}

\paragraph{Proof.}  The Erd\H{o}s-Turan inequality (Proposition ~\ref{pr32})
gives that for any $K$, 
\beql{453}
{\Bbb E}_{\omega_{N}}[ D(\{ y_{1},...,y_{N}\})] \le
\frac{1}{K+1} +
 3 \sum_{k=1}^{K}\frac{1}{kN} 
{\Bbb E}_{\omega_{N}}[ |\hat{U}_{N} (k, \omega)| ].
\eeq
By the Cauchy-Schwarz inequality, \eqn{431}, and the definition of the
two-dimensional Diophantine 
property  we have that 
\beql{454}
{\Bbb E}_{\omega_{N}} [ |\hat{U}_{N}(k, \omega)|]
\le \Big({\Bbb E}_{\omega_{N}}\Big[ |\hat{U}_{N}(k, \omega)|^{2}\Big]
\Big)^{\frac 12} \le (1+C(\theta_1,\theta_2)^{-2} 
k^{2\alpha})^{\frac 12}  \sqrt{N}.
\eeq
Using this in \eqn{453} we obtain that 
for an appropriate constant $C_1(\theta_1,\theta_2)$,
$$
E_{\omega_{N}}[ D(\{ y_{1},...,y_{N}\})] \le
\frac{1}{K+1} + C_1(\theta_1,\theta_2) \frac{K^{\alpha}}{\sqrt{N}}.
$$
Choosing $K=N^{\frac{1}{2(1+\alpha)}}$ we obtain the Theorem. 
$~~~\bsq$ \\

\paragraph{Remark.} The stochastic process studied in
this section can be reformulated in terms of the iterates of a
skew-product dynamical system, as defined in
 Cornfeld, Fomin and Sinai \cite[Chap. 10]{CFS82}
and Petersen \cite{Pe89}.
  Let $\Sigma = \{0,1\}^{{\NN}}$
denote the set of all zero-one sequences
$\bs = (s_{0}, s_{1}, s_{2}, ...)$, with the product
topology, which is a compact space with natural
invariant measure, and let  $S: \Sigma \to \Sigma$ be
the shift operator 
$S(s_{0}, s_{{1}}, s_{2}, ...)= (s_{1}, s_{2}, s_{3},...)$.
The skew-product dynamical system
 $T: \Sigma \times \TT \to \Sigma \times \TT$ 
over the base $\Sigma$, with fibers $\TT = \RR/\ZZ$, is defined  by
 $$
 T(\bs, x) := (S(\bs), x+ f(s_{0}) (\bmod 1)),
 $$
  with
 $f(0)= \theta_{1}, f(1)= \theta_{2}$, respectively. Here the
 initial condition is $(\bs^{(0)}, x_{0})$, 
with $\bs^{(0)} \in \Sigma$ being a random
 starting point.  The invariant measure on
$\Sigma \times \TT$ is the product measure, using 
 Lebesgue measure on $\TT$,
 and $T$ is ergodic with respect to this measure 
if at least one of $\theta_{1}$
 and $\theta_{2}$ is  irrational.
 The initial result of this section (Theorem ~\ref{th41})
shows weak convergence of almost all orbits to Lebesgue measure
on $\TT$ for the dynamical system. This result is true
in great generality for ergodic skew products.
However the detailed result on rate of
convergence to Lebesgue measure (Theorem ~\ref{th42}) 
relies on specific properties
of this dynamical system.

%***************************************************
%
%
% Section 5 Application to 3X+1 Problem
%
%
%***************************************************

\section{Application to the $3x+1$ map}

We can describe the $3x+1$ iteration applied to an
integer $m$ in terms of the parity of its iterates.
We set $T^{(0)} (m) =m$ and define the {\em parity sequence} 
$\{b_{k}(m): k \ge 0\}$ with 
each $b_k(m) \in \{ 0, 1\}$ by
\beql{501}
b_{k} (m) \equiv  T^{(k)}(m)~(\bmod~2).
\eeq

%*******************************************
% proposition 5.1
%*******************************************
%

\begin{prop}~\label{pr51}
(1) The $k$-th iterate $T^{{(k)}}(m)$ for $k \ge 1$ has the form
\beql{502}
T^{(k)}(m) = \frac{3^{b_{0}(m)+... + b_{k-1}(m)}}{2^{k}}m + R_{k}(m)
\eeq
in which the remainder term
\beql{503}
R_{k}(m) := \sum_{j=0}^{k-1} b_{j}(m) 
\frac{3^{b_{j+1}(m)+...+ b_{k-1}(m)}}{2^{k-j}}
\eeq
depends only on $m \pmod{2^{k}}$.

(2) Each $b_{k}(m)$ depends only on $m \pmod{2^{k+1}}$.
For each  vector $(b_{0}, b_{1}, ..., b_{{N-1}}) \in \{0, 1\}^{N}$
there is a unique residue class $m ~(\bmod~2^{N})$ such that
\beql{504}
b_{k}(m) = b_{k}  ~~\mbox{for}~~ 0 \le k \le N-1.
\eeq
\end{prop}

\paragraph{Proof.}
(1) This is easily proved by  induction on $k$, see 
Lagarias \cite[(2.6)]{La85}.

(2) This is also proved by induction on
$k$, see  Lagarias \cite[Theorem B]{La85}. $~~~\bsq$ \\

We define $x_{k}(m) = T^{(k)}(m)$ 
and view
\beql{506}
\tilde{x}_{k}(m) := \frac{3^{b_{0}(m)+... + b_{k-1}(m)}}{2^{k}}m
\eeq
as an approximation to $x_{k}(m)$. 
Viewing the base  $B\ge 2$ as fixed, we set
$y_{k}(m) := \log_{B} x_{k}(m) $
and the main result will concern the discrepancy of most
sets $\sY_{N}(m) := \{ y_1(m), ..., y_{N}(m) \}$. 
We approximate the $y_k(m)$ by 
\begin{eqnarray}\label{507a}
\tilde{y}_{k}(m) & := &\log_{B} \tilde{x}_{k}(m) 
 =  \log_{B} m + \Big(\sum_{j=0}^{k-1} b_{j}(m)\Big) \log_{B} 3
 -  k \log_{B} 2.
\end{eqnarray}
and we will  study the
sets 
 $\tilde{\sY}_{N}(m) := \{ \tilde{y}_{1}(m), ..., \tilde{y}_{N}(m) \})$
for variable $m$ as realizations of a stochastic process of
the kind treated in \S4.

The following lemma shows that  the error of approximation of
$\sY_{N}(m)$ by $\tilde{\sY}_{N}(m)$
is exponentially small in $N$ for most $m$.

%*******************************************
% lemma 5.1
%*******************************************
%
\begin{lemma}~\label{le51}  Let the integer $B \ge 2$ be fixed.   
There exists  an 
exceptional subset $E_B(N)$ of integers $1\le m \le 2^N$ such that 
$$
|E_B(N)| \le 2^{1+\frac{99}{100}N}, 
$$ 
and such that if $1\le m\le 2^N$ is not in $E_B(N)$ then 
\beql{511}
|y_{k}(n) - \tilde{y}_{k}(n)| \le 2^{1-\frac{1}{100}N}
~~\mbox{for}~~ 1\le k \le N,
\eeq
for every $n\equiv m \pmod {2^N}$.
\end{lemma}

\paragraph{Proof.} We will prove more, and show that the set 
$E_B(N)$ may be taken to be the set of integers $1\le m \le 2^N$ 
such that either $m\le 2^{\frac{99}{100}N}$, or 
$b_0(m)+\ldots+b_{N-1}(m) \le \frac 25 N$.  Since all $2^N$ possible 
choices for the parities $b_0(m)$, $\ldots$, $b_{N-1}(m)$ 
occur exactly once, 
we see that the number of $m$ satisfying the second criterion above is 
$\le \sum_{j\le \frac 25 N} {\binom{N}{j}} \le 2^{H(\frac{2}{5})N}
\le 2^{\frac{99}{100}N}$, where 
$H(x)=-x \log_2 x -(1-x)\log_2(1-x)$ is the binary entropy function.
Thus $|E_B(N)| \le 2^{1+\frac{99}{100}N}$, as desired.  It remains now to 
show \eqn{511} holds for $m\notin E_B(N)$.  

Suppose now that $m\notin E_B(N)$ and that $n\equiv m \pmod{2^N}$.  
Proposition ~\ref{pr51} gives that $b_{k-1}(n)= b_{k-1}(m)$ 
and $R_{k}(n)= R_{k}(m)$
for $1 \le k \le N$. Observe that
$$
\frac{x_k(n)}{\tilde{x}_k(n)} 
= 1 + \frac{R_{k}(n)}{\tilde{x}_k(n)} = 1 +\frac{R_k(m)}{\tilde{x}_k(n)}  
\le 1+ \frac{R_k(m)}{\tilde{x}_k(m)} = \frac{x_k(m)}{\tilde{x}_k(m)}, 
$$
from which it follows that $y_k(n) -\tilde{y}_k(n) \le 
y_k(m) - \tilde{y}_k(m)$. 
Thus it suffices to verify \eqn{511} for $n=m$.  

 From \eqn{503} we see that 
$$
R_k(m) \le \sum_{j=0}^{k-1} \frac{3^{k-j-1}}{2^{k-j}} \le 
\Big(\frac 32\Big)^k.
$$ 
Applying this bound together with $\log (1+\xi) \le \xi$, 
we obtain that  
\beql{515}
y_{k}(m) - \tilde{y}_{k}(m) = \log_B 
\Big( 1 +\frac{R_k(m)}{\tilde{x}_k(m)}\Big) 
\le  \frac{1}{\log B} \frac{R_k(m)}{\tilde{x}_k(m)} 
\le \frac{1}{\log B} \frac{1}{m} 3^{k-b_0(m)-\ldots-b_{k-1}(m)}.
\eeq
Since $m\notin E_B(N)$ we have that $m >2^{\frac{99}{100}N}$, and 
in addition that 
$$
k-\sum_{j=0}^{k-1} b_j(m) =\sum_{j=0}^{k-1} (1-b_j(m)) \le 
\sum_{j=0}^{N-1} (1-b_j(m)) \le N- \frac{2}{5} N = \frac{3}{5} N.
$$
Thus from \eqn{515} we deduce for $m\notin E_B(N)$ that 
$$
y_k(m)-\tilde{y}_k(m) 
\le \frac{1}{\log B} 2^{-\frac{99}{100}N} 3^{\frac 35 N} \le 
2^{1- \frac 1{100}N}, 
$$
(since   $3^{\frac {3}{5}}< 2^{\frac{98}{100}}$) 
which proves the Lemma.  $~~~\bsq$ \\

We wish to bound the discrepancy of 
most sets $\tilde{\sY}_{N}(m)$, viewed over a range $1 \le m \le X$, with
$X \ge 2^{N}$.  We will study the translated sets 
\beql{520}
\tilde{\sY}_{N}^{\ast}(m) :=\tilde{\sY}_{N}(m) - \log_{B} m,
\eeq
so that the initial element $\tilde{y}_0^{\ast}(m)$ is zero.  
Since the discrepancy function is translation invariant 
we have that
\beql{520b}
D(\tilde{\sY}_{N}^{\ast}(m) )= D( \tilde{\sY}_{N}(m) )
\eeq
Note also that $\tilde{\sY}_N^{\ast}(m) = \tilde{\sY}_N^{\ast}(m+2^N)$ 
and so it will suffice to consider the range $1 \le m \le 2^{N}$.

%*******************************************
% lemma 5.2
%*******************************************
%

\begin{lemma}~\label{le52}
Let $B \ge 2$ and $N \ge 1$ be fixed. Then the ensemble  
$
\{\tilde{\sY}_{N}^{\ast}(m) : ~1 \le m \le 2^{N} \}
$
of $2^N$ sequences of length $N$
is identical in distribution
with the distribution $\omega_{N}$ of
the first $N$ elements of the stochastic process 
$\sP(\theta_1,\theta_2, y_{0}=0)$,
with parameters
$\theta_1 = \log_{B} \frac{3}{2}$  and
 $\theta_2= \log_{B} \frac 12$.
\end{lemma}

\paragraph{Proof.} From the definitions we see easily that 
$\tilde{y}_k^{\ast}(m) = \tilde{y}_{k-1}^{\ast}(m) + \theta_1$ 
if $b_{k-1}(m) =1$, and that $\tilde{y}_k^{\ast}(m) = 
\tilde{y}_{k-1}^{\ast}(m) + \theta_2$ 
if $b_{k-1}(m)=0$. 
Proposition ~\ref{pr51}(2) shows that 
for $1 \le m \le 2^{N}$ all possible patterns $(b_{0}, b_{1}, ..., b_{N-1})$
occur exactly once. This corresponds exactly to independent draws in the
stochastic process $\sP(\theta_1, \theta_2, y_{0}=0)$; the $2^{N}$ possible
sequences $\omega_{N}$ of length $N$ of  $\sP(\theta_1, \theta_2, y_{0}=0)$
have equal probabilities and match the sequences above.
$~~~\bsq$ \\

%We now show that for any real $B > 1$ the pair 
%$(\theta_{1}, \theta_{2})= (\log_{B}\frac{3}{2}, \log_{B}\frac 12)$ 
%are two-dimensional Diophantine.
%*************************************************************************
%
%    Lemma 5.3 dioph approx properties
%
%*************************************************************************
\begin{lemma} ~\label{le53}  For each real $B > 1$
the  pair $(\theta_1,\theta_2)=(\log_B \frac 32, \log_B\frac 12)$ 
is two-dimensional Diophantine with exponent $7.616$.
\end{lemma}
\paragraph{Proof.}  We invoke a result of Rhin \cite{Rh87} 
(see inequality (8) there) 
obtained using Pad{\' e} approximation methods:  There exists a positive 
constant $C$ such that for integers $u_0$, $u_1$ and $u_2$ with 
$\max(|u_1|,|u_2|) \ge 1$ 
we have 
\beql{523}
|u_0 + u_1 \log 2 + u_2 \log 3| \ge C  \Big( \max(|u_1|,|u_2|)\Big)^{-7.616}.
\eeq 

Let $k$ be a large positive integer and suppose that $\ell_1$ is 
the nearest integer to $k\theta_1$ and that $-\ell_2$ is 
the nearest integer to $k\theta_2$.  Thus $|k\theta_1-\ell_1| = 
\Vert k\theta_1 \Vert$ and $|k\log_B 2 -\ell_2| = |k\theta_2 
+\ell_2|=\Vert k\theta_2 \Vert$.  Note that both $\ell_1$ and 
$\ell_2$ are positive and roughly of size $k$.  On the one hand we 
see that 
$$
\Big| \frac{\log (3/2)}{\log 2} 
-\frac{\ell_1}{\ell_2} \Big| 
= \Big|  \frac{\ell_2 k \log_B (3/2) - \ell_1 k \log_B 2}{k \ell_2 \log_B 2} 
\Big| \le 
\frac{\ell_2 \Vert k\theta_1 \Vert + 
\ell_1 \Vert k\theta_2 \Vert}{k\ell_2 \log_B 2}.
$$
On the other hand we see that by \eqn{523}
$$
\Big| \frac{\log(3/2)}{\log 2} -\frac{\ell_1}{\ell_2}\Big| 
= \Big| \frac{\ell_2 \log 3 - (\ell_1+\ell_2) \log 2}{\ell_2 \log 2}\Big| 
\ge C \frac{(\ell_1 +\ell_2)^{-7.616}}{\ell_2 \log 2}.
$$
Since $\ell_1$ and $\ell_2$ are roughly of 
size $k$, combining the above two statements immediately gives the 
Lemma. $~~~\bsq$

% line 1.1482

%*******************************************
% Proof of Theorem 2.1
%*******************************************
%

\paragraph{Proof of Theorem \ref{th21}.}  We view 
the integer $B \ge 2$ and $N \ge 1$ as fixed.  
Consider the realizations $\omega_N$ of the stochastic process 
$\sP(\theta_1,\theta_2,y_0=0)$ 
with $\theta_1=\log_B \frac 32$ and $\theta_2 = \log_B \frac 12$.  
By Lemma~\ref{le53} and Theorem~\ref{th42} 
we obtain that (with $\alpha = 7.616$) 
$$
{\Bbb E}_{\omega_N} [ D(\{y_1,\ldots,y_N\}) ] \le 
C N^{-\frac{1}{2(1+\alpha)}} \le CN^{-\frac{1}{18}},
$$
for an appropriate positive constant $C$.  Using
Markov's inequality that
${\rm Prob}[Y \ge a] \le \frac{{\Bbb E}[Y]}{a}$
for a nonnegative random variable $Y$, we  deduce that 
$$
{\rm Prob}\Big[ D(\{y_1,\ldots,y_n\}) \ge N^{-\frac{1}{36}} \Big] 
\le C N^{-\frac{1}{36}}.
$$
Invoking Lemma~\ref{le52} we conclude  
that the exceptional set of $m$ with $1 \le m \le 2^{N}$ such that
$D(\tilde{\sY}_{N}(m)) \ge    N^{- \frac{1}{36}}$
has cardinality at most $C  N^{- \frac{1}{36}} 2^{N}$.
By Lemma~\ref{le51} we know that for most $1\le m\le 2^N$ 
the sets $\sY_{N}(m)$ and $\tilde{\sY}_N(m)$ are very close 
term by term, and by  Proposition~\ref{pr33} for such $m$ the 
discrepancies $D(\sY_N(m))$ and $D(\tilde{\sY}_N(m))$ are 
very nearly equal.
Thus we may deduce that the exceptional set of 
$m$ with $1 \le m \le 2^{N}$ such
that
$$
D(\sY_{N}(m)) \ge    N^{- \frac{1}{36}} + 2^{2- \frac{1}{100} N}
$$
has cardinality at most
$$
C N^{- \frac{1}{36}} 2^{N} +  2^{1+\frac{99}{100}N}.
$$
This easily gives the conclusion of the theorem for $X=2^{N}$.

It remains to treat the case $X > 2^{N}$.  Suppose 
$\ell 2^{N}< X \le (\ell+1)2^{N}$, 
for some $\ell \ge 1$.  Since the discrepancies $D(\tilde{\sY}_{N}(m))$ 
are periodic $\pmod{2^N}$ we see that the exceptional set of 
$m\le X$ with large discrepancy contains no more than $\ell+1$ times 
the number of exceptional $m\le 2^N$.  
This completes the proof.  $~~~\bsq$ \\

{\bf Acknowledgements.}  We thank Jens Marklof for  
interesting discussions on this problem, and  the reviewer
for helpful comments and references.
 Both authors are partially supported by
NSF grants.

%\clearpage

\noindent Jeffrey C. Lagarias \\
Dept. of Mathematics\\
The University of Michigan \\
Ann Arbor, MI 48109-1043\\
\noindent email: {\tt lagarias@umich.edu}\\

\noindent K. Soundararajan  \\
Dept. of Mathematics\\
The University of Michigan \\
Ann Arbor, MI 48109-1043\\
\noindent email: {\tt ksound@umich.edu}


\begin{thebibliography}{99}
\bibitem{Be38}
F. Benford,
The law of anomalous numbers,
Proc. Amer. Phil. Soc. {\bf 78} (1938), 551--572.

\bibitem{Be05}
A. Berger, 
Multi-dimensional dynamical systems and Benford's law,
Discrete Contin. Dyn. Sys. {\bf 13} (2005), 219--237.

\bibitem{BBH05}
A. Berger, L. Bunimovich and T. Hill,
One-dimensional dynamical systems and Benford's law,
Trans. Amer. Math. Soc. {\bf 357} (2005), 197--219.

\bibitem{BP00}
K. Borovkov and D. Pfeifer, Estimates for the Syracuse problem
via a probabilistic model,
Theory of Probability and its Applications {\bf 45}, No. 2 (2000),
300--310.

\bibitem{CFS82}
I. P. Cornfeld, S. V. Fomin and Ya. G. Sinai,
{\it Ergodic Theory}, 
Springer-Verlag, New York 1982.

\bibitem{Di77}
P. Diaconis,
The distribution of leading digits and uniform distribution
mod 1,
Ann. Prob. {\bf 5} (1977), 72--81.


\bibitem{DT97}
M. Drmota and R. F. Tichy,
{\em Sequences, Discrepancies and Applications},
Lecture Notes in Math. 1651, Springer-Verlag: New York 1997.

\bibitem{He79}
M-R. Herman,
Sur la conjugaison diff\'{e}rentiable des diff\'{e}omorphismses 
du cercle \`{a} des
rotations, Publ. Math. IHES, No. 49 (1979), 5--233.


\bibitem{Hi95}
T. Hill,
Base-invariance implies Benford's law,
Proc. Amer. Math. Soc. {\bf 123} (1995), 187--195.

\bibitem{Hi96}
T. Hill,
A statistical derivation of the significant digit law,
Statistical Science {\bf 86} (1996), 358--363.

\bibitem{KM04}
A. V. Kontorovich and S. J. Miller,
Benford's law, values of $L$-functions and the $3X+1$ problem,
Acta Arithmetica {\bf 120} (2005), 269--297.
%(eprint: {\tt arXiv math.NT/0412003})

\bibitem{KS03}
A. V. Kontorovich and Ya. G. Sinai,
Structure theorem for $(d,g,h)$-maps,
Bull. Braz. Math. Soc. {\bf 33} (2002), 213--224.

\bibitem{KN74}
L. Kuipers and H. Niederreiter,
{\it Uniform Distribution of Sequences},
John Wiley \& Sons: New York 1974.

\bibitem{La85}
J. C. Lagarias,
The $3x+1$ problem and its generalizations,
Amer. Math. Monthly {\bf 92} (1985), 3--23.

\bibitem{LW92}
J. C. Lagarias and A. Weiss,
The $3X+1$ problem: two stochastic models,
Ann. Applied Prob. {\bf 2} (1992), 329--361.

\bibitem{Ma03}
J. Marklof,
Pair correlation densities of inhomogeneous quadratic forms,
Ann. Math. {\bf 158} (2003), 419--471.

\bibitem{Ma02}
J. Marklof,
Pair correlation densities of inhomogeneous quadratic forms II,
Duke Math. J. {\bf 115} (2002) 409--434. 
(Correction, ibid {\bf 120}
(2003), 227--228).

\bibitem{Ma05}
J. Marklof,
Mean square value of exponential sums related to the representation
of integers as sums of squares,
Acta Arith. {\bf 117} (2005), no. 4, 353--370.

\bibitem{Mo94}
H. L. Montgomery,
{\em Ten Lectures on the Interface Between Analytic Number Theory
and Harmonic Analysis},
CBMS Regional Conference Series No. 84,
Amer. Math. Soc.: Providence, RI 1994.

\bibitem{New81}
Simon Newcomb,
Note on the frequency of use of the different digits in natural numbers,
Amer. J. Math. {\bf 4} (1881) 39--40.

\bibitem{Pe89}
K. R. Petersen,
{\it Ergodic Theory}, Corrected reprint,
Cambridge Univ. Press: Cambridge 1989.

\bibitem{Rh87}
G. Rhin,
Approximants de Pad\'{e} et mesures effectives d'irrationalit\'{e},
in: {\it S\'{e}minaire de Th\'{e}orie des Nombres, Paris 1985-1986},
pp. 155--164,
Progress in Mathematics, Vol. 71,
Birkh\"{a}user, Boston 1987. 






\bibitem{Ro04}
E. Roosendaal, private communication ($n= 72,19136, 41637,72362,71195
\approx 7.2 \cdot 10^{21}$).

\bibitem{Sch71} 
W. Schmidt, Approximation by algebraic numbers,
Enseignement Math. {\bf 17} (1971), 187--253.

\bibitem{Si03}
Ya. G. Sinai,
Statistical $(3x+1)$ problem. Dedicated to the memory of
J\"{u}rgen K. Moser, 
Comm. Pure Appl. Math. {\bf 56} (2003), 1016--1028.




\bibitem{Wi98}
G. J. Wirsching,
{\it The Dynamical System Generated by the $3n+1$ Function},
Lecture Notes in Math. 1681, Springer-Verlag: New York 1998.

\bibitem{Yo84}
J.-C. Yoccoz, 
Conjugaison diff\'{e}rentiable des diff\'{e}omorphismes du 
cercle dont le nombre de
rotation v\'{e}rifie une condition diophantienne, 
Ann. Sci. \'{E}cole Norm. Sup.
 (4) {\bf 17} (1984), 333--359.
 
 \bibitem{Yo02}
 J.-C. Yoccoz,
 Analytic linearization of circle diffeomorphisms,
 in: {\it Dynamical systems and small divisors (Certraro, 1998))},
 Lecture Notes in Math. 1784, 
 Springer-Verlag: New York 2002, pp. 125--173.
\end{thebibliography}
\end{document}